\documentclass[11pt]{amsart}
\usepackage{graphicx}
\usepackage{amsmath,amssymb}
\newtheorem{thm}{Theorem}[section]
\newtheorem{cor}[thm]{Corollary}
\newtheorem{lem}[thm]{Lemma}

\theoremstyle{definition}

\theoremstyle{remark}
\newtheorem{rem}[thm]{Remark}
\theoremstyle{Example}

\numberwithin{equation}{section}

\begin{document}
\title{A new family of  elliptic curves with positive ranks arising from the Heron triangles}
\author[F.A. Izadi   ]{F.A. Izadi   }
\address{ Mathematics Department Azarbaijan university of  Tarbiat Moallem ,
 Tabriz, Iran  f.izadi@utoronto.ca  farzali.izadi@gmail.com}

\author[$\thickspace $ F. Khoshnam ]{$\quad $F. Khoshnam }
\address{ Mathematics Department Azarbaijan university of  Tarbiat Moallem ,
 Tabriz, Iran khoshnam@azaruniv.edu}
\author[$\thickspace $ K. Nabardi ]{K. Nabardi }
\address{ Mathematics Department Azarbaijan university of  Tarbiat Moallem ,
 Tabriz, Iran nabardi@azaruniv.edu}

\begin{abstract}
The aim of this paper is to introduce a new family of elliptic curves with positive ranks.
These elliptic curves have been constructed with certain rational numbers, namely a, b, and
c as sides of Heron triangles having rational areas $k$. It turned out  that the torsion
groups of this family are of the form $\frac{\Bbb{Z}}{2\Bbb{Z}}\times \frac{\Bbb{Z}}{2\Bbb{Z}}$
and also the rank is positive.
\end{abstract}
\maketitle {\small{\bf Keywords:} Elliptic curves, Ranks, Torsion
group, Heron's Triangle.}

\section{Introduction}
As is well-known, an elliptic curve E over a field $\Bbb{K}$ can
be explicitly expressed by the generalized Weierstrass equation of the form $$E
: y^2+a_1xy+a_3y=x^3+a_2x^2+a_4x+a_6$$ where,
$a_1,a_2,a_3,a_4,a_6\in \Bbb{K}$.  In this paper we are interested
in the case of
$\Bbb{K}=\Bbb{Q}$.\\
By  the Mordell-Weil theorem \cite{sil},  every elliptic curve over
$\Bbb{Q}$ has a commutative group $E(\Bbb{Q})$ which is finitely
generated, i.e., $E(\Bbb{Q})\cong\Bbb{Z}^r\oplus
E(\Bbb{Q})_{tors} $, where $r$ is a nonnegative integer and
$E(\Bbb{Q})_{tors}$ is the subgroup of elements of finite order
in $E(\Bbb{Q})$. This subgroup is called the torsion subgroup of
$E(\Bbb{Q})$ and the integer $r$ is called the rank of E
and is denoted by the rank(E).\\
By Mazur's theorem,  the torsion subgroup $E(\Bbb{Q})_{tors}$ is
one of the following 15 forms: $\mathbb{Z}/n\mathbb{Z}$ with
$1\leq n\leq 10$ or $n=12$, $\mathbb{Z}/2\mathbb{Z}\times
\mathbb{Z}/2m\mathbb{Z}$ with $1\leq m\leq 4$. Besides, it is not
known which values of rank $r$ are possible. The folklore
conjecture is that a rank can be arbitrarily  large, but it seems to
be very hard to find examples with large ranks. The current record
is an example of elliptic curve over $\mathbb{Q}$ with rank$\geq$
28, found by Elkies  in May 2006 (see \cite{el}). Having
classified the torsion part, one interested in seeing whether or
not the rank is unbounded among all the elliptic curves. There is
no known guaranteed algorithm to determine the rank and it is not
known which integers can occur as ranks.
\section{NEW FAMILY OF ELLIPTIC CURVES}
In our work, we are going to study the rank of a new family of
elliptic curves based on some triangles with rational areas $k$
and rational sides. Recall that any triangle with rational sides
and rational area is simply called a Heron triangle due to Heron
of Alexandria (C. 10 A.D. -C. 75 A.D.) and  that Heron's formula
for the area, say $S$,  of a triangle with sides $(a,b,c)$ is
given by
$$S=\sqrt{P(P-a)(P-b)(P-c)}$$
Where, $P=\frac{(a+b+c)}{2}$ is  the semi perimeter.\\
Take $S=v$ and $P=u$, so we have $v^2=u(u-a)(u-b)(u-c)$,  where $a$, $b$, and $c$ are the
sides of a Heron triangle.\\
With the change of coordinates, $(u,v)\rightarrow
(\frac{1}{\zeta},\frac{\eta}{{\zeta}^2})$, we obtain $${\eta}^2=(1-a\zeta)(1-b\zeta)(1-c\zeta)$$
with one more change of coordinates, $(\zeta,\eta)\rightarrow
(\frac{-x}{abc},\frac{y}{{abc}})$, we get

$$y^2=(x+ab)(x+bc)(x+ac).$$ We will show that
$$E(\Bbb{Q})_{tors}=\frac{\Bbb{Z} }{2\Bbb{Z}}\oplus \frac{\Bbb{Z}
}{2\Bbb{Z}}=\{ (-ab,0),(-bc,0),(-ac,0),\infty\}.$$ However from
$x=0$  we get $y=\pm abc$. We will show that this two
points $(0,\pm abc)$ are belong to the free torsion part of
$E(\Bbb{Q})$.

Now we define the elliptic curve associated to the triple
\{a(k),b(k),c(k)\} arising from a Heron triangle with rational area $k$ and sides

\[  \left\lbrace
  \begin{array}{c  l}
   a(k)=5k^2-4k+4,& \\
                &\\
    b(k)=\frac{1}{2}k(k^2-4k+20), & \\
                 &\\
    c(k)=\frac{1}{2}(k+2)(k^2-4), &
  \end{array}
\right. \] given by N. J. Fine \cite{fine},

i.e.,
$$y^2=(x+a(k)b(k))(x+b(k)c(k))(x+a(k)c(k)).$$
Then the curve E has three rational points of order two :
\[  \left\lbrace
  \begin{array}{c  l}
    P1=[\frac{1}{2}(5k^2-4k+4)k(k^2-4k+20),0],&      \\
                                              &\\
    P2=[\frac{1}{2}(5k^2-4k+4)(k+2)(k^2-4),0],& \\
                                            &\\
     P3=[\frac{1}{2}k(k^2-4k+20)(k+2)(k^2-4),0]. & \\

     &
  \end{array}
\right. \] We can make a change of coordinates, $(x,y)\rightarrow
(x-a(k)c(k),y)$ which sends $(0,0)\rightarrow
(a(k)c(k),0)$. Obviously such a change won't affect the structure
of the group $E(\mathbb{Q})$. Thus, given the restriction that we
are considering curves with a rational two-torsion point, we can
assume that $E$ is given by $y^2=x^3+Ax^2+Bx$. This leads to the
following values of the coefficients $A$ and $B$:
$$A=\frac{1}{4}k^6-16k^4-3k^5+96k^3-44k^2-16k+32,$$
\begin{equation}
\begin{split}
  B &=
  -\frac{1}{4}(k^6-4k^4-12k^5+96k^3-16k^2-192k+64)\\ &\times
   (15k^4-72k^3+40k^2-32k-16)
   .
\end{split}
\end{equation}

The corresponding value of the discriminant is
\begin{equation}
\begin{split}
  \Delta &=
  16(A^2-4B)B^2=\frac{1}{16}k^2(k^2-4k+20)^2(k^3-8k^2+4k-16)^2\\ &\times
   (k^2-12k+4)^2(k-2)^4(k+2)^4(5k^2-4k+4)^2(3k^2-12k-4)^2
   .
\end{split}
\end{equation}
From the discriminant of the curve we see that for the values
$k=0$, $2$,  and -2 curve reduces to singular form.

\begin{thm}
Let $a(k)$, $b(k)$ and $c(k)$ be defined as above. then the
elliptic curve $$E :
y^2=\Big(x+a(k)b(k)\Big)\Big(x+b(k)c(k)\Big)\Big(x+a(k)c(k)\Big)$$
over $\mathbb{Q}(k)$ has the torsion group isomorphic to
$\frac{\Bbb{Z}}{2\Bbb{Z}}\times \frac{\Bbb{Z}}{2\Bbb{Z}} $.
\end{thm}
\begin{proof}
With coordinate transformation $x\longrightarrow x-a(k)c(k)$ our
curve leads to
$$y^2=x\Big(x-a(k)b(k)+a(k)c(k)\Big)\Big(x-b(k)c(k)+a(k)c(k)\Big).$$
 Kubert \cite{4} showed that if $y^2=x(x+r)(x+s)$, with $r,s\neq0$ and $s\neq
r$, then the torsion subgroup is $\frac{Z}{2Z}\times\frac{Z}{2Z}$.
Therefore for our curve it is sufficient to prove that
$-a(k)b(k)+a(k)c(k)\neq 0$, $-b(k)c(k)+a(k)c(k)\neq 0$, and
$-a(k)b(k)+a(k)c(k)\neq -b(k)c(k)+a(k)c(k)$ for all $k$ in
$\mathbb{Q}(k)$. Firstly we try to find all  $k$ in
$\mathbb{Q}(k)$ that satisfy in the equalities, $a(k)=b(k)$,
$a(k)=c(k)$,  and $b(k)=c(k)$. Theses equalities give rise to the following equalities:
\[  \left\lbrace
  \begin{array}{c  l}
    k^3-14k^2+28k-8=0,&      \\
                                              &\\
    k^3-8k^2+4k+16=0,& \\
                                            &\\
     6k^2-24k-8=0. & \\

     &
  \end{array}
\right. \] By the Eisenstein   theorem if $k=\frac{p}{q}\in
\mathbb{Q}$ for relatively primes $p$ and $q$, be a solution to the first equation,
the possibilities for $p$ are $\{\pm 1,\pm 2,\pm 4,\pm 8\}$ which
shows that  $k=2$ satisfies  this cubic equation, but we have already assume that  $k\neq 0,2,-2.$\\
Also for the second equation, we see that there is no solution.\\
For the third equation it's trivial that the solutions are $k=2\pm
\frac{4}{3}\sqrt{3}$.\\
Therefore, for every $k$ in $\mathbb{Q}$, $a(k)\neq b(k)$,
$c(k)\neq b(k)$, $c(k)\neq a(k)$ and the torsion group for main
curve is $\frac{\Bbb{Z}}{2\Bbb{Z}}\times \frac{\Bbb{Z}}{2\Bbb{Z}}
$.
\end{proof}

\begin{cor}
$rank(E(\mathbb{Q}))\geq 1.$
\end{cor}
\begin{proof}
By Theorem $2.1$, the point
$$P_k=(0,abc)=(0,1/4(5k^2-4k+4)k(k^2-4k+20)(k+2)(k^2-4))$$
on $E(\mathbb{Q})$ has infinite order, which shows that
$rank(E(\mathbb{Q}))\geq 1$.
\end{proof}


%
\begin{lem}
For each $1\leq r \leq 6$ there exists a $k$ such that the
elliptic curve $$y^2=(x+ab)(x+bc)(x+ac)$$ has a torsion group
isomorphic to $\frac{\Bbb{Z}}{2\Bbb{Z}}\times
\frac{\Bbb{Z}}{2\Bbb{Z}} $ and its rank equals to $r$.
\end{lem}
\begin{proof}
It is trivial by the results given in the table above.
\end{proof}
\begin{rem}
If $a$, $b$ and $c$ are the sides of right triangle then, the relations
$a^2+b^2=c^2$, $a^2+c^2=b^2$ and $b^2+c^2=a^2$ lead to $k=6$, $2/3$ that has the curves with ranks equal 1. It is
remarkable that in our computation we only meet this two curves with rank =1.
\end{rem}
\section{Our Computation}
In this stage we want to find curves having large ranks possible. The main idea here is that a curve is more
likely to have large rank if $| E(\mathbb{Q})|$ is relatively
large for many primes $p$. We will use the following realization
of this idea. For a prime $p$ we put $a_p =a_p(E) = p + 1 - |
E(\mathbb{F}_p)|$ and
$$S(N, E) =\sum_{p\leq N, p     \medspace  prime}(1-\frac{p-1}{ |
E(\mathbb{F}_p)|})\log (p)=\sum_{p\leq N, p\medspace
prime}(\frac{-a_p+2}{p+1-a_p})\log (p).$$ It is experimentally
known that one may expect the high rank curves have large $S(N,
E)$. Some arguments show that the Birch and Swinnerton-Dyer
conjecture gives support to this observation. The sum $S(N, E)$
can be very efficiently computed (e.g., using PARI )for $N <
10000$. After this sieving method, we may continue to investigate
the best, let us say, $1$ percent of the curves. Since, we are working
with curves with torsion points of order 2, we may compute the Selmer rank which is
the best known upper bound for the actual rank of the curves. In our
computations, we used the sage software \cite{sage} and Cremona's
MWRANK program \cite{CER} for computing the Mordell-Weil rank of
the curves. Following is a table of elliptic curves with
corresponding values of a, b, c and
 ranks =1, 2, 3, 4, 5 and 6 . Also we have calculated the rank of this family for $100000$ curves
 in the range of  $1\leq k \leq 100 $ using the MWRANK. The results appear in the following table.
\begin{table}
[!h]\begin{center}
\begin{tabular}{|c|c|c|} \hline
$k$ & Curve & Rank \\
\hline \hline
& & \\
& $y^2=x^3+\frac{-3859986810117979136}{59604644775390625}x^2+\frac{-302381696902314690275394654830592}{1136868377216160297393798828125}x  $&\\
$\frac{98}{625}$&&6\\ 
&$ +\frac{720840680923373992917188523670698174399532498944}{21684043449710088680149056017398834228515625}$&\\ 
&&\\ \hline
& & \\
$11$&\footnotesize $y^2=x^3+4547347x^2+\tiny{6818384095380}x+3363133863125708100  $&5\\
&&\\\hline 
& & \\
$19$&\footnotesize $y^2=x^3+89515187x^2+\tiny{2523432520031220}x+22674567869155130588100  $&4\\
&&\\ \hline 
& & \\
$3$&\footnotesize $y^2=x^3+6899x^2+\tiny{14152500}x+8901922500  $&3\\
&&\\ \hline 
& & \\
$4$&\footnotesize $y^2=x^3+26432x^2+\tiny{225607680}x+613652889600  $&2\\
&&\\ \hline 
& & \\
$6$&\footnotesize $y^2=x^3+192512x^2+\tiny{12079595520}x+247390116249600  $&1\\
&&\\\hline 

\hline
\end{tabular}
\end{center}
A table of  elliptic curves  with ranks = 1,2,3,4,5,6 .
\end{table}


\begin{table}
\begin{center}
\begin{tabular}{|c|c|} \hline
Rank & Percent  \\
\hline
$6$&$0.5\%$\\
\hline
$5$&$3.7\%$\\
\hline
$4$&$23.4\%$\\
\hline
$3$&$36.7\%$\\
\hline
$2$&$16.2\%$\\
\hline
undetermined &$19.5\%$\\
\hline

\end{tabular}
\end{center}
The table for  percents.
\end{table}


\newpage

\begin{thebibliography}{99}
\bibitem {CER} {\sc J.Cremona,}, {\it mwrank program},
http://maths.nottingham.ac.uk/personal/jec/ftp/progs/.
\bibitem{el} {\sc N. D. Elkies}, {\it $\mathbb{Z}^{28}$ in $E(\mathbb{Q})$}, ect., Number Theory Listserver, May 2006.
   \bibitem{fine} {\sc N. J. Fine}, {\it on rational triangles}, American Mathematical Monthly. 83, no. 7(1976)
     517-521, MR0414484, zbI 0341.10016.

\bibitem{4} D. S. Kubert, Universal Bounds On The Torsion Of
Elliptic Curves, Proc. London Math.Soc. (3), 33, 1976,pp.193-237.
    \bibitem {sage} {\sc Sage software}, {\it Version 4.3.5}, http://sagemath.org .
\bibitem {sil} {\sc  J. H. Silverman}, {\it Advanced Topics in the Arithmetic
of Elliptic Curves}, Springer-Verlag, New York, 1994.



\end{thebibliography}
\end{document}